\documentclass[twoside]{article}
\usepackage{amssymb}
\input{diagrams}
\usepackage{theorem}

\makeatletter
\@addtoreset{equation}{subsection}
\@addtoreset{equation}{section}
\makeatother
\newtheorem{theorem}[equation]{Theorem}
\newtheorem{proposition}[equation]{Proposition}
\newtheorem{corollary}[equation]{Corollary}
\newtheorem{lemma}[equation]{Lemma}
{\theorembodyfont{\rmfamily}
\newtheorem{statement}[equation]{}
\newtheorem{definition}[equation]{Definition}

\newtheorem{remark}[equation]{Remark}
\newtheorem{example}[equation]{Example}
}

\renewcommand
        {\theequation}
        {\arabic{section}.\ifnum\value{subsection}=0\else\arabic{subsection}.\fi\arabic{equation}}

\newcommand{\Ol}{{\cal O}}

\newcommand{\Z}{{\mathbb Z}}

\newcommand{\R}{{\mathbb R}}
\newcommand{\op}{{\cal O}_{\mathbb P^1}}

\newcommand{\opn}{{\cal O}_{\mathbb P^n}}
\newcommand{\oqn}{{\cal O}_{\mathbb Q^n}}
\newcommand{\f}{\varphi}
\newcommand{\ra}{\rightarrow}
\newcommand{\lra}{\longrightarrow}

\newcommand{\E}{{\cal E}}

\newcommand{\Proj}{{\mathbb P}}
\newcommand{\Segno}

\date{}
\title{Ample vector bundles with sections vanishing on special varieties}
\author{Marco Andreatta, Gianluca Occhetta}
\begin{document}
\maketitle 
\begin{abstract}
Let $\cal E$ be an ample vector bundle of
rank $r$ on a complex projective manifold $X$ such that there exists a section $s \in \Gamma(\cal E)$
whose zero locus $Z = (s = 0)$ is a smooth submanifold of the expected dimension $\textrm{dim~} X - r: =
n -r$. Assume that $Z$ is not minimal; we investigate the hypothesis under which the extremal contractions of $Z$ can be
lifted to $X$. Finally we study in detail the cases in which $Z$ is a scroll, a quadric bundle or a del Pezzo fibration.
\end{abstract}
\begin{figure}[bp]
\footnotesize{1991 {\it Mathematics Subject Classification}: Primary 14E30, 14E40; Secondary 14F05 14J45. }
\end{figure}
\setlength{\parindent}{0pt}

\section{Introduction}

A very classical and natural way of classifying complex projective manifolds $X$
consists in slicing $X$ with a number of general hyperplane sections obtaining
in this way a complex manifold of smaller dimension which is likely classifyable.
Then one should {\sl ascend} the geometrical properties of this new manifold
and obtain a complete description of $X$. To stress the classical flavor of
this approach it is sometime called {\sl Apollonius method} \cite{Fub}.
The hard part of the Apollonius method are the ascending properties; in this
paper we will consider this problem in a slightly more general set up.\par\vspace{0.4 cm}

Let $\cal E$ be an ample vector bundle of
rank $r$ on $X$ such that there exists a section $s \in \Gamma(\cal E)$
whose zero locus $Z = (s=0)$ is a smooth submanifold of the expected dimension $dim X - r: = n -r$.\par

Assume that $Z$ is not minimal in the sense of Mori's theory,
that is $-K_Z$ is not nef; thus $Z$ has at least one extremal ray \cite[Cone Theorem]{Mo1} and an
associated extremal contraction \cite[Kawamata-Shokurov base point free theorem]{KMM}.

Our question will then be, under which condition this contraction
can be {\sl lifted} to the ambient variety, determining its
structure; this general situation is studied in section $3$;
suppose that $F_Z$ is an extremal face in $\overline{NE(Z)}$ with
supporting divisor $K_Z +\tau H_Z$; a lifting property is proved
under the assumption that
\begin{equation} \label{ampleH}
H_Z \quad \textrm{is the restriction of an ample line bundle~} H \textrm{~on~} X .
\end{equation}
Next we discuss some special situations in which the assumption (\ref{ampleH})
can be avoided.
In the rest of the paper we consider some special cases, namely if $Z$
is a scroll, a quadric bundle or a del Pezzo fibration; the results are described in
theorems (\ref{mainscroll}), (\ref{mainquadric}), (\ref{maindelpezzo}) and corollaries.

These results generalize classical ones by L. B\v{a}descu (see
\cite{Ba1}, \cite{Ba2}, \cite{Ba3}) and A.J. Sommese (see \cite{So1} and chapter 5 of \cite{BS2}, in particular
Theorem 5.2.1 which should be compared with our results in section 3) and
more recent ones by A. Lanteri and H. Maeda
(\cite{LM1},\cite{LM2},\cite{LM3}), who were the first ones to study the problem of special sections of ample
vector bundles. \\
We would like to thank J. Wi{\'s}niewski for some useful remarks and
the referee, who pointed out some inaccuracies in the first version of the manuscript.
During the preparation of this paper we were partially supported by  the MURST of the Italian Government.

\section{Notations and generalities}

We use the standard notation from algebraic geometry.
In particular it is compatible with that of \cite{KMM}
to which we refer constantly; we suggest to the reader also the survey
\cite{AW}.
We just explain some special definitions
and propositions used frequently.

In this paper $X$ will always stand for a smooth complex projective variety
of dimension $n$. Let $Div(X)$ the group of Cartier divisors on $X$;
denote by $K_X$ the {\sf canonical divisor} of $X$, an element of
$Div(X)$ such that $\Ol_{X}(K_X) = \Omega^n_{X}$.
Let $N_1(X)=\frac{\{1-cycles\}}{\equiv}\otimes \R$,
$N^1(X)= \frac{\{divisors\}}{\equiv}\otimes \R$ and
$\overline {<NE(X)>}=\overline{\{\mbox{effective 1-cycles}\}}$;
the last is a closed cone in $N_1(X)$. Let also $\rho(X)=dim_{\R}N^1(X)$.

Suppose that $K_X$ is not nef, that is there exists an effective curve $C$
such that $K_X\cdot C<0$.

\begin{theorem}\cite{KMM}
Let $X$ as above and $H$ a nef Cartier divisor such that
$F:= H^{\bot} \cap \overline {<NE(X)>} \setminus \{0\}$
is entirely contained in the set
$\{Z\in N_1(X) :K_X\cdot Z<0\}$,
where $H^{\bot} = \{Z:H\cdot Z=0\}$.
Then there exists a projective morphism $\f:X\ra W$ from $X$ onto a normal
variety $W$ which is characterized by the following properties:
\begin{itemize}
\item[{i})]  For an irreducible curve $C$ in $X$, $\f(C)$
is a point if and only if $H.C = 0$, if and only if
$cl(C) \in F$.
\item[{ii})] $\f$ has only connected fibers
\item[{iii})] $H = \f^*(A)$ for some ample divisor $A$
on $W$.
\end{itemize}
\label{contractionth}
\end{theorem}

\begin{definition} (\cite{KMM},
definition 3-2-3). Using the notation of the above theorem,
the map $\f$ is called a {\sf Fano-Mori contraction}
(or an
{\sf extremal contraction}); the set $F$ is an {\sf extremal face},
while the Cartier
divisor $H$ is a {\sf supporting divisor} for the map $\f$ (or the face $F$).
The contraction is of fiber type if
$dim W < dim X$, otherwise it is birational.
If $dim_{\R}F = 1$ the face $F$ is called an {\sf extremal ray}, while $\f$
is called an {\sf elementary contraction}.
\end{definition}

\begin{remark}\label{nefvalue}
Note that a supporting divisor for a Fano-Mori contraction is of
the form $H=K_X +rL$ where $r$ is a positive integer. In fact if $H$
is a supporting divisor then $H-K_X$ is an ample line bundle by
Kleiman's criterion.
\end{remark}

\begin{remark}\label{biraz} Let $\pi:X\ra V$ denote a contraction of an
extremal
face $F$, supported by $H=\pi^*A$. Let
 $R$ be an extremal ray in $F$ and $\rho:X\ra W$ the contraction of $R$.
Then $\pi$ factors trough $\rho$ (this is because $\pi^*A\cdot R=0$).
\end{remark}

\begin{remark} We have also \cite{Mo1} that if $X$ has an extremal ray $R$
then there exists a rational curve $C$ on $X$ such that
$0< -K_X \cdot C\leq n+1$ and
$R=R[C]:=\{D\in <NE(X)>: D\equiv \lambda C, \lambda\in \R^+\}$.
Such a curve is called an {\sf extremal curve}.
\end{remark}

The last remark was generalized by P. Ionescu and J. Wi\'sniewski as in the
follow.

\begin{definition} Let $\f$ be a Fano-Mori contraction
of $X$ and let $E = E(\f)$ be the exceptional locus of $\f$
(if $\f$ is of fiber type then $E:=X$);
let $S$ be an irreducible component of a (non trivial) fiber $F$.
We define the positive integer $l$ as
$$l =  min \{ -K_X\cdot C: C \textrm{~is a rational curve in~} S\}.$$
If $\f$ is the contraction of a ray $R$, then $l$ is called the length of
the ray.
\end{definition}

\begin{proposition}\label{iowi} \cite{Wi1} In the set-up of the previous
definition the
following formula holds
$$dim S + dim E \geq dimX + l -1.$$
\label{diswis}
\end{proposition}

In particular this implies that if
$\varphi$ is of fiber type then $l  \leq (dim Z- dim W +1)$ and if
$\varphi$ is birational then $l \leq (dim Z -dim (\varphi(E))$.\par
\medskip

If a manifold admits different extremal contractions, then the dimensions of different general fibers are bounded
by the following

\begin{theorem}(\cite[Theorem 2.2]{Wi1})\label{manyrays} Let a manifold $X$ of
dimension $n$ admit k different contractions (of different
extremal rays). If by $m_i, i=1,2,\dots,k$ we denote dimensions of
images of these contractions, then
$$ \sum^k_{i=1} (n-m_i)\le n$$
\end{theorem}

\medskip
It is very useful to understand when a contraction is elementary;
for this we will use in this paper the following result:

\begin{proposition}\cite{ABW1} \label{elementary} Let $\pi:X\ra W$ be a
contraction of
a face such that $dim X > dim W = m$. Suppose that for every rational curve $C$
in a general fiber of $\pi$ we have
$-K_X\cdot C\geq (n+1)/2$.
Then $\pi$ is an elementary contraction except if
\begin{itemize}
\item[a)] $-K_X\cdot C=(n+2)/2$ for some rational curve $C$ on $X$,
$W$ is a point, $X$ is a Fano manifold
of pseudoindex $(n+2)/2$ and $\rho(X)=2$
\item[b)]  $-K_X\cdot C=(n+1)/2$
for some rational curve $C$, and $dim$W$\leq 1$
\end{itemize}
\label{fibelementare}
\end{proposition}

{\bf Proof} \quad For the reader's convenience we will give a sketch of the proof. \\
Let $T$ be a non breaking family of rational curves of $X$, which is dominant and whose dimension at every point is 
greater or equal to $-K_X.l -2$, where $l$ is a curve in $T$.
The existence of this family follows as in \cite{Mo2} from the fact that through a general point of a general fiber
there passes a rational curve $l$ such that 
$$ -K_X.l = -K_F .l \le n-m+1$$
and that, on the other hand, by our assumptions, $-K_X .l > (n-m+1)/2$.\\
Suppose that $\pi$ is not an elementary contraction. Therefore there is a contraction
$\varphi = \textrm{cont}_R$ where $R=\mathbb R_+[l_2]$ is an extremal ray contracted by $\pi$
and not containing $l$.
Let $F$ be a fiber of $\varphi$; then \quad $\textrm{dim}~F \ge l(R) -1$\quad (\ref{iowi}); let $T_p$ be the locus of curves from the 
family $T$ which pass through a given point $p$; we have dim~$T_p \ge -K_X.l -1$.
But, by the non-breaking lemma (see \cite{Wi2}), we must have
$$ \textrm{dim}~F + \textrm{dim}~T_p \le \textrm{dim}~X $$
if $T_p \cap F \not = \emptyset$ and for $p \notin F$, that is
$$ -K_X.l -K_X.l_2 \le n+2$$
If $-K_X.C > (n+2)/2$ for every rational curve in a general fiber we have a contradiction.\\
If $-K_X.C = (n+2)/2$ then we have equality everywhere; in particular we have that
dim~$F= -K_X.l -1 = n -\textrm{dim~}T -1$. Then, by \cite[lemma 1.4.5]{BSW} we have that
$NE(X) = NE(F) + \mathbb R_+[l]$. Since $F$ is a positive dimensional fiber of an elementary
contraction, we conclude that $NE(F) = \mathbb R_+$ and thus that $\rho(X) = 2$. Thus 
$\rho(Y) = 0$, i.e. $Y$ is a point; therefore $X$ is a Fano manifold of pseudoindex $(n+2)/2$.\\
If $-K_X.C = (n+1)/2$ then a general fiber of $\pi$, $G$, is a Fano manifold of pseudoindex
$\ge (n+1)/2$. If $m \ge 2$ we have that the pseudoindex of $G$ is greater or equal to $(\textrm{dim}~G)/2 +1$;
therefore $\rho(G) = 1$ (see \cite{Wi2}).
This implies in particular that, if $\pi$ is not an elementary contraction, then $\varphi = \textrm{cont}_R$, as above,
is birational and therefore that \quad dim~$F \ge -K_X.l_2$\quad (\ref{iowi}); since $\rho(X) \ge 3$ it follows as above that
$$ -K_X.l -K_X.l_2 > n,$$
and thus we arrive at the contradiction with $-K_X.C \le (n+1)/2$ using again \cite[lemma 1.4.5]{BSW}.

\begin{definition} Let $L$ be an
ample line bundle on $X$. The pair $(X,L)$ is called a scroll (respectively
a quadric fibration, respectively a del Pezzo fibration)
over a normal variety
$Y$ of dimension $m < n$ if there exists a surjective morphism
with connected fibers
$\phi: X \ra Y$ such that
$$K_X+(n-m+1)L \approx p^*{\mathcal L}$$
(respectively $K_X+(n-m)L \approx p^*{\mathcal L}$,
respectively $K_X+(n-m-1)L \approx p^*{\mathcal L}$)
for some ample line bundle ${\mathcal L}$ on $Y$.
$X$ is called a classical scroll or a $\mathbb P$-bundle (respectively
quadric bundle) over a projective
variety $Y$ of dimension $r$ if there exists a surjective morphism
$\phi : X\ra Y$ such that every fiber is isomorphic to $\Proj^{n-r}$
(respectively to a quadric in $\Proj^{n-r+1}$)
and if there exists a vector bundle $\cal E$ of rank $n-r+1$ (respectively
of rank $n-r+2$) on $Y$ such that  $X\simeq \Proj(\cal E)$
(respectively exists an embedding of $X$ over $Y$ as a divisor of $\Proj(\cal E)$ of relative degree 2).
\end{definition}

\begin{remark} A scroll is a Fano Mori contraction of fiber type
such that the inequality in (\ref{diswis}) is actually an equality,
i.e. $l  = (dim X- dim Y +1)$ and moreover if $C$ is a rational curve
such that $-K_X^.C = l$ then it exists an ample line bundle $L$ such that
$L^.C = 1$, i.e. $C$ is a line with the respect to $L$.
The contrary is almost true in the sense that if $\f$ is a Fano Mori
contraction
with the above properties then it factors through a scroll, that is the
face which
is contracted by $\f$ contains a sub-face whose contraction is a scroll.

Similarly a quadric (respectively a del Pezzo) fibration is a
Fano Mori contraction of fiber type
such that $l  = (dim X- dim Y)$ (resp.  $l = (dim X- dim Y-1)$)
and moreover if $C$ is a rational curve
such that $-K_X^.C = l$ then it exists an ample line bundle $L$ such that
$L^.C = 1$, i.e. $C$ is a line with the respect to $L$.
\end{remark}

\begin{theorem}\cite[Proposition 3.2.1]{BS2}\label{pbundle}
Let $p:X \to Y$ be a surjective equidimensional morphism onto a normal variety $Y$ and let
$L$ be an ample line bundle on $X$ such that $(F,L_F)
\simeq(\mathbb P^r, {\cal O}(1))$ for the general fiber $F$ of
$p$. Then $p:X \to Y$ gives to $(X,L)$ the structure of a $\mathbb
P^d$-bundle.
\end{theorem}

\begin{remark} Let $\f: X \ra Y$ be a scroll
and let $\Sigma \subset Y$ be the set of points $y$ such that $dim(\f^{-1}(y))
> k:= dim X- dim Y$ then $codim \Sigma \geq 3$, thus if $dim Y \leq 2$ then
the scroll is a $\Proj^k$-bundle
(\cite[Theorem 3.3]{So}).
\end{remark}

\begin{theorem}\cite[Theorem B]{ABW3}\label{qbundle} Let $(X,L)$ be a quadric
fibration $\f: X \to Y$, with $L$ an ample line bundle on $X$;
assume that $\f$ is an elementary contraction and that $\f$ is
equidimensional. Then ${\cal E}:= p_*L$ is a locally free sheaf of
rank dim $X -$ dim $Y +2$ and $L$ embeds $X$ into $\mathbb P({\cal
E})$ as a divisor of relative degree $2$, i.e. $X$ is a classical
quadric bundle.
\end{theorem}

\begin{remark}
Let $\varphi:X \to Y$ be a scroll (respectively a quadric fibration, respectively a del Pezzo fibration) and let dim~$X =  n$, dim~$Y = m$; it follows
directly from (\ref{elementary}) that if $n\ge 2m -1$ (respectively if $n \ge 2m +1$ and $(m,n) \not =  (0,2),(1,3)$, respectively
if $n\ge 2m+3$ and $(m,n) \not = (0,4),(1,5)$) then $\varphi$ is an elementary contraction, i.e. the contraction of an extremal ray.
\end{remark}

\begin{lemma}\label{det} Let $\cal E$ be an ample vector bundle of rank $r$
on a complex variety $X$. For any rational curve $C \subset X$
we have
$$ (det {\cal E}).C \ge r.$$
Moreover, if $C$ is smooth and the equality holds, then ${\cal E}_C = {\cal
O}_{\mathbb P^1}(1)^{\oplus r}$.
\end{lemma}

\begin{lemma}\label{edim} Let $Y$ be a complex projective variety of
dimension $n$, $\cal E$ an
ample vector bundle on $Y$, $s$ a global section of $\cal E$;
denote with $V(s)$ the zero set of $s$. Then
$$\textrm{dim}V(s)\ge n-r .$$
\end{lemma}
{\bf Proof.} \quad See \cite[Example 12.1.3]{Ful}.

\begin{proposition}\label{intersezione} Let $X$, $\cal E$ and $Z$ be as before.
Let $Y$ be a subvariety of $X$ of
dimension $\ge r$. Then dim~$Z \cap Y \ge$~dim~$Y -r$ .\end{proposition}

{\bf Proof.} Consider ${\cal E}_Y$, the restriction of $\cal E$ to $Y$ and
$s_Y$, the
restriction of $s$ to $\Gamma(Y,{\cal E}_Y)$. Applying lemma
(\ref{edim}) to $Y$ and $s_Y$ we get
$$ \textrm{dim~}(Z \cap Y) = \textrm{dim~}V(s_Y)\ge \textrm{dim~}Y -r $$

\section{Lifting of contractions}

\begin{statement} \label{setup}
Let $\cal E$ be an ample vector bundle of
rank $r$ on $X$ such that there exists a section $s \in \Gamma(\cal E)$ whose zero locus
$Z = (s=0)$ is a smooth submanifold of the expected dimension $dim X - r: =
n -r$. Note that, with this assumptions, the restriction of $\cal E$ to $Z$ is the normal bundle
$N_{X}Z$ \cite[Example 6.3.4]{Ful}.
\end{statement}

The idea of this section is to investigate the relation between $\overline{NE(X)}$
and $\overline{NE(Z)}$; one basic result in this direction is the following Lefschetz type theorem
proved by Sommese in \cite{So2} and with slightly weaker assumptions in \cite{LM1}.

\begin{theorem}\label{lef}
Let $X$, $\cal E$ and $Z$ be as in \ref{setup} and let $i:Z~\hookrightarrow~X$
be the embedding. Then\par\vspace{0.3 cm}
(\ref{lef}.1) $H^i(i):H^i(X,\mathbb Z) \to H^i(Z, \mathbb Z)$ is an isomorphism
for $i \le \textrm{dim}~Z -1$\par
(\ref{lef}.2) $H^i(i)$ is injective and its cokernel is torsion
free for $i= \textrm{dim~}Z$\par
(\ref{lef}.3) $\textrm{Pic}(i): \textrm{Pic}(X)\rightarrow \textrm{Pic}(Z)$
is an isomorphism for $dim Z \ge 3$.\par
(\ref{lef}.4) $\textrm{Pic}(i)$ is injective
and its cokernel is torsion free for $dim Z = 2$.\par
(\ref{lef}.5) $\rho(X) = \rho (Z)$ for $dim Z \ge 3$.\par
\end{theorem}

Note that, although the Picard groups are isomorphic, in general the ample
cone of $X$ is
properly contained in the ample cone of $Z$. However, in special cases,
something can be said;
the following proposition generalizes a result of J. A. Wi\'sniewski on
divisors.\par

\begin{proposition}\label{fano} Let $X$ be a Fano manifold of dimension $n$, $\cal E$
an ample vector bundle
of rank $r$ on $X$ and $Z$ the zero locus of a section of $\cal E$, smooth
and of the expected
dimension. If $Pic(X) \cong Pic(Z)$ and $X$ has no elementary extremal
contractions
with all fiber of dimension $\leq r$
then a line bundle on $X$ is ample if and only if its restriction to $Z$ is
ample.
The assumption is satisfied for instance if
all the extremal rays of $X$ have length
$ l(R) \ge r+2 $ or $ l(R) \ge r+1$ if $R$ is not nef.
\par
\end{proposition}

{\bf Proof.} \quad Observe that, since $X$ is Fano, a line bundle $\cal L$
on $X$
is ample if and only if it has positive intersection with any extremal ray
of $X$.\par
So take a line bundle ${\cal L}_Z$ ample on $Z$; if we prove that every
extremal ray of $X$
contains the class of a curve lying on $Z$ we can conclude that $\cal L$ is
ample on
$X$.\par
By our assumption for every extremal ray of $X$ its associated contraction
has a fiber $F$ of
dimension $\geq  r+1$; thus
$$ \textrm{dim} F + \textrm{dim~} Z \ge n+1$$
and therefore, by proposition~(\ref{intersezione}) the intersection of $Z$ and $F$
contains a curve, which belongs to the ray $R$.
Using (\ref{iowi}) one shows immediately that the assumption on the length
implies
the lower bound on the fiber.$\square$\par

\vspace{0.3cm}

The same idea allows us to prove the following\par

\begin{theorem}\label{lifting}(Lifting of contractions) Let $X$, $\cal
E$ and $Z$ be as in (\ref{setup}) and assume that $Z$ is not
minimal. Let $F_Z$ be an extremal face of $Z$ and $D_Z=K_Z+\tau
H_Z$ a good supporting divisor of $F_Z$. Assume that there exists
an ample line bundle $H$ on $X$ which is the extension of $H_Z$.
Then $D = K_X + det E + \tau H$ is nef, but not ample;
thus it defines an extremal face $F_X$ of $X$. Moreover, if $\tau
\ge 2$ and $dim Z \ge 3$, under the identification of $N_1(X)$
with $N_1(Z)$ we have $F_X = F_Z$ and the contraction of every ray
spanning $F_Z$ lifts.
\end{theorem}

{\bf Proof}. \quad Suppose that $D$ is not nef.
There exists a curve $C$ on $X$ such that $D.C < 0$; therefore
there exists an extremal ray $R$ on $X$ on which $D$ is negative
and s.t. $l(R) \ge r+ \tau+ 1$.\par

From the inequality (\ref{iowi})

$$ dim F \ge l(R) - 1 $$
where $F$ is a non trivial general fiber of the contraction of $R$ and
this yields
$$ dim F + dim Z \ge r + \tau + n - r = \tau + n \geq n+1$$
so, recalling that we can assume $\tau \ge 1$ (\ref{nefvalue}),
in view of proposition~(\ref{intersezione}) a curve of the ray $R$
lays on $Z$ and this is absurd, since $D_{|Z}$ is nef.\par To prove
the last claim, observe that every extremal ray $R$ in the face
$F_X$ has length $l(R) \ge \tau +r$, so the general non trivial
fiber of the contraction of $R$ has dimension dim~$F \ge \tau +r-1
\ge r+1$, so, in view of proposition (\ref{intersezione}) a curve of the
ray $R$ lays on $Z$. $\square$

\begin{proposition} \label{nohypo}
The hypothesis on the ampleness of $H$ is not necessary if $dim Z
\geq 2$ and $Pic(Z) \cong \mathbb Z$ or, more generally, if
$Pic(i):Pic(X)\to Pic(Z)$ is an isomorphism and $\overline{NE(Z)}=
\overline{NE(Z)}_{K_Z\ge 0} + R$, i.e. if $Z$ has only one extremal
ray. Example (\ref{ex1}) shows that the hypothesis is necessary if
$Z$ has at least two extremal rays.
\end{proposition}

{\bf Proof.} \quad There exists a line bundle $L$ which is ample on
$X$; the restriction of this line bundle to $Z$, $L_Z$ is ample on
$Z$, so, if $K_Z$ is not nef there exist a rational number $\sigma
>0$ such that $K_Z + \sigma L_Z$ is nef but not ample and it
defines an extremal face $G_Z$ (\cite[Kawamata rationality theorem]{KMM}).
But we are supposing that on $Z$
there is only one extremal ray, thus $F_Z = G_Z$.\par

\begin{remark}\label{nocor} Note that the proof of (\ref{nohypo}) actually shows there is always an
extremal contraction on $Z$ which can be lifted to $X$.\end{remark}

\begin{lemma} \label{nohyposp}
If $\varphi:Z \to W$ is a $\mathbb P$-bundle contraction 
on a smooth minimal variety $W$
then $Z$ has only one extremal ray.
\end{lemma}

{\bf Proof.} \quad Suppose that $Z$ has another extremal ray,
$R_1$; there exists a rational curve $C_0$ such that $-K_Z.C_0
> 0$ and $\f(C_0)$ is not a point.
Let $C=\f(C_0)$, let $\nu:\mathbb P^1 \to C$ be the normalization
of $C$ and consider the fiber product

\begin{equation}
\begin{diagram}
Z \times_W \mathbb P^1& \rTo^{\bar\nu} & Z \\
\dTo^{\bar\varphi~} & & \dTo_{~\varphi} \\
\mathbb P^1 & \rTo^{\nu}& W\\
\end{diagram}
\end{equation}
\vspace{0.3 cm}

$\bar\f: Z_C:= Z \times_W \mathbb P^1 \to \mathbb P^1$ is a
$\mathbb P$-bundle on $\mathbb P^1$ and so $\rho(Z_C)=2$; the morphism $\bar\nu$ induces a map of
spaces of cycles $N_1(Z_C) \to N_1(Z)$ which is an embedding. The
Mori cone $NE(Z_C)$ is contained in the intersection $N_1(Z_C) \cap NE(Z)$ and
so, since $N_1(Z_C)$ is a plane in $N_1(Z)$ and passes through two
different extremal rays of $Z$, $NE(Z_C)$ is contained in the
negative part of $NE(Z)$.\par By \cite[Corollary 2.8]{KoMiMo}
$-K_{Z_C/\mathbb P^1}$ is not ample, so there exist an horizontal
curve $C_1$ on $Z_C$ such that $-K_{Z_C/\mathbb P^1}.C_1 \le 0$;
noting that
$$-K_{Z_C/\mathbb P^1} = \bar\nu^*\f^*K_W - \bar\nu^* K_Z$$
we get
$$\bar\nu^*\f^*K_W.C_1 = -K_{Z_C/\mathbb P^1}.C_1+\bar\nu^* K_Z.C_1 \le K_Z.\bar\nu(C_1) < 0$$
and therefore $K_W.\f(\bar\nu(C_1)) < 0$, which contradicts the
minimality of $W$.\par

\begin{remark}  We found the idea of the proof of (\ref{nohyposp}) in \cite{SW2} and \cite{KoMiMo}.\end{remark}

\vspace{0.3cm}
For the rest of this section we will be in the hypothesis of
theorem (\ref{lifting}) and we will denote by $\f: Z \ra W$
the contraction of the face $F_Z$ and by $\phi : X \ra Y$ the contraction of
$F_X$.
Let also $m= dim W$.

By the adjunction formula
 $-K_Z = -(K_X+\textrm{det} {\cal E})_Z$, so $-K_Z$ is $\phi$-ample.\par

On $Z$ we have thus two contractions, $\varphi$ and $\phi_Z$. Now we
are going to investigate the relation between them.
Clearly we have a commutative diagram\par

\begin{equation}\label{diagram}
\begin{diagram}
X & \lTo^i & Z \\
\dTo^{\phi~} &\ldTo^{\phi_Z} & \dTo_{~\varphi} \\
Y& \lTo^{\pi}& W\\
\end{diagram}
\end{equation}
\vspace{0.3 cm}

\begin{lemma}\label{surjectivity} $\phi_Z(Z) \supseteq \phi(E(\phi))$.
\end{lemma}

{\bf Proof.}\quad We reason as in the proof of (\ref{lifting}):
since $\phi$ is the contraction of a ray of length
$l(R)\ge r +\tau$,  a non trivial fiber of $\phi$ has dimension
$\geq \tau +r -1$ and thus it has nonempty intersection with $Z$.
\par

\begin{proposition}
\label{fiber} If the contraction $\varphi$ is of fiber type
then also $\phi$ is of fiber type.\end{proposition}

{\bf Proof.} \quad If $\varphi$ is of fiber type the commutativity of the
diagram
(\ref{diagram}) implies that also $\phi_Z$ is of fiber type, so $Z$
is contained in the exceptional locus of $\phi$, $E(\phi)$, and
by lemma (\ref{surjectivity}) $\phi_Z(Z)=\phi(E(\phi))$.\par
Suppose that $\phi$ is birational; in this case $E(\phi) \subsetneq X$
and
$$ \textrm{dim~}\phi(E(\phi)) = \textrm{dim~}\phi_Z(Z) < \textrm{dim~} Z =
n-r.$$
$Y$ has dimension $n$, so it is possible to find a subvariety $Y' \subset Y$
of dimension $r$ which has empty intersection with
$\phi((E(\phi))$; away from $E(\phi)$, $\phi$ is an isomorphism, so
$X' = \phi^{-1}(Y') \subset X$ is a subvariety of $X$ of dimension
$r$ which has empty intersection with $E(\phi)$ and therefore with
$Z$, but this is absurd by proposition (\ref{intersezione}). \quad $\square$
\par

\begin{proposition}\label{agreement} If $\f$ is of fiber type
or if $\tau \ge 2$, $\phi_Z$ has connected fibers.
Moreover $\phi_Z$ factors as $\phi_Z= \sigma \circ \varphi$ where $\sigma: W \to \phi_Z(Z)$ is the
normalization morphism. In particular, if $\varphi$ is of fiber type then $\phi_Z = \varphi$.
\end{proposition}

{\bf Proof.} \quad The fibers of $\phi_Z$ are of the form $Z \cap F$
with $F$ fiber of $\phi$.
If $\varphi$ is of fiber type, then the same is for $\phi$ (see \ref{fiber})
whose fibers have thus dimension $\geq dim X- dim Y = n-m$;
so dim~$Z \cap F\ge n - r - m\ge 1$. If $\tau \ge 2$, reasoning as in the
proof of theorem (\ref{lifting}) we again have dim~$Z \cap F\ge 1$.
So theorem (\ref{lef}.1)
applies to $F$ and ${\cal E}_F$ and gives $H^0(Z \cap F,\mathbb Z) \cong
H^0(F,\mathbb Z)\cong
\mathbb Z$. Using the Universal Coefficient Theorem we get $H_0(Z\cap F)
\cong \mathbb
Z$.\par
Let $\sigma: W' \to \phi_Z(Z)$ be the normalization of $\phi_Z(Z)$; by the universal property
$\phi_Z$ factors through $\widetilde{\phi_Z}:Z \to W'$ and $\sigma$; note also that, since
$\phi_Z$ has connected fibers, the same is true for $\widetilde{\phi_Z}$.

Let $C \subset Z$ be any curve contracted by $\phi _Z$ (and hence by $\widetilde{\phi_Z}$); thus
$(K_X +det{\cal E} +\tau H)_Z.C= 0$ which is equivalent to
$(K_Z + \tau H_Z).C = 0$, i.e. $C$ is contracted by $\f$.
By the commutativity of the diagram every curve contracted by $\f$ is
contracted by $\phi _Z$ (and hence by $\widetilde{\phi_Z}$).
Therefore $\varphi$ and $\widetilde{\phi_Z}$ are two Fano Mori contractions which contract the 
same extremal face, so they are the same morphism. \\
To prove the last claim recall that, by lemma (\ref{surjectivity}) and proposition (\ref{fiber}) if 
$\varphi$ is of fiber type then $\phi_Z(Z) = Y$ and hence $\phi_Z(Z)$ is normal. \par

\begin{proposition}\label{birubiru}
If the contraction $\varphi$ is birational and $\tau \ge 2$, then
also $\phi$ is birational.
\end{proposition}

{\bf Proof.} \quad Suppose $\phi$ is of fiber type; reasoning again as
in the proof of (\ref{lifting}) we can prove that dim~$Z\cap F \ge 1$
for the generic fiber of $\phi$, so $\phi_Z=\varphi$ is of fiber type, a
contradiction in view of proposition (\ref{agreement}).\par

\begin{remark} If $\varphi$ is birational and $\tau = 1$ then $\phi$ can be of fiber type (see case 3. of proposition (\ref{F1}).
\end{remark}

\section{Scrolls and $\mathbb P^d$-bundles}

\begin{theorem}\label{mainscroll}
Let $X$, $\cal E$ and $Z$ be as in (\ref{setup}) with $dim Z \geq
2$. We assume that $Z$ has a scroll contraction $\f: Z \ra W$ with
respect to an ample line bundle on $Z$, $H_Z$, which is the
restriction of an ample line bundle $H$ on $X$.

Then $X$ has a Fano-Mori contraction $\phi:X\to W$ which is of fiber type and
with supporting divisor $D=K_X+det {\cal E}+(n-m-r+1)H$.
The general fiber of $\phi$ is isomorphic to $\mathbb P^{n-m}$ and ${\cal E}$ restricted to
it is $\oplus^r{\cal O}_{\mathbb P}(1)$.

If $\varphi$ is elementary or $dim X = n \geq 2 m-1 = 2 dim W -1$
(this is always the case if $dim W \leq 3$)
then $\phi$ is elementary and it is a scroll contraction (i.e. it
is supported by the divisor $K_X+(n-m+1)H$) and moreover in the
second case even $\f$ had to be elementary.
\end{theorem}

{\bf Proof.}\quad The morphism $\varphi$ is a contraction supported
by $K_Z +(n-r-m+1)H_Z$, so, applying theorem (\ref{lifting}),
we get a contraction $\phi:X\to Y$, defined by an high multiple of
$D=K_X+det {\cal E}+ (n-m-r+1)H$; this contraction is of fiber type
and $Y = W$ by proposition (\ref{agreement}). Let $F$ be a general
fiber of $\phi$; then $F$ is a smooth Fano manifold of dimension
$n-m$ such that $-K_F=({det \cal E}+(n-m-r+1)H)_F$. Thus $F=\mathbb
P^{n-m}$ and $\cal E$ restricted to it is $\oplus^r{\cal
O}_{\mathbb P}(1) (see \cite[]{Pe1}).$ Moreover, for any line in a general fiber $(det
{\cal E} -r H).l = 0$.\par
 Assume now that dim~$X \ge 2$ dim~$W -1$; note that dim~$X \ge$ dim~$Z+1
\ge$ dim~$W+2$,
so the inequality holds for $dim W \leq 3$.\par By the proposition
(\ref{elementary}) the contraction $\phi:X\to W$ is an elementary
contraction and so
$$det {\cal E}= rH + \phi^*B$$
that is $\phi$ is supported by $K_X +(n-m+1)H$; note that also
$\varphi:Z \to W$ had to be elementary, by the last claim of
theorem (\ref{lifting}) if dim $Z \ge 3$ and by
(\ref{elementary}) if dim $Z =2$. $\square$\par

\begin{corollary}\label{mainbundle}
Assume now that $Z = \Proj ({\cal F})$ for some vector bundle
${\cal F}$ on $W$, and its tautological bundle is the restriction
of an ample line bundle on $X$; then $X =
\Proj ({\cal G})$ for some vector bundle ${\cal G}$ on $W$ which
admits ${\cal F}$ as a quotient; in this case ${\cal E}= \xi_{\cal G} \otimes \phi^*{\cal I}{\check{~}}$
where ${\cal I}$ fits into the exact sequence
$$ 0 \lra {\cal I} \lra {\cal G} \lra {\cal F} \lra 0$$
and $\phi: X \to W$ is the $\mathbb P$-bundle contraction.
\end{corollary}

{\bf Proof.} \quad The theorem gives us a contraction $\phi:X \to
W$; we claim that $\phi$ is equidimensional; in fact if it has any
fiber of dimension $> n-m$ then, by proposition (\ref{intersezione}),
even $Z
\ra W$ should have a fiber of dimension $>(n-m-r)$. Since $\f: Z
\ra W$ is elementary $\phi$ is a scroll with the respect to $H$. The
first part of the corollary is proven by (\ref{pbundle}). The second part is a well known
fact about vector bundles (see \cite[B.5.6.]{Ful}).

\begin{example} Let ${\cal F}$ be an ample vector bundle on a smooth curve $C$ of genus $g >0$. If ${\cal F}$ is
decomposable into a sum of $r$ bundles ${\cal F}_i$, by \cite[Corollary 4.20]{Fu2} each ${\cal F}_i$
fits into an exact sequence
$$ 0 \lra {\cal O}_C \lra {\cal G}_i \lra {\cal F}_i \lra 0 $$
with ${\cal G}_i$ ample; so we can construct an exact sequence
$$ 0 \lra \oplus^r{\cal O}_C \lra {\cal G} = \mathop{\oplus}\limits^r_{i=1} {\cal G}_i \lra {\cal F} = \mathop{\oplus}\limits^r_{i=1} {\cal F}_i  \lra 0.$$
On $X = \mathbb P({\cal G})$ the vector bundle ${\cal E} = \xi_{\cal G} \otimes p^*(\oplus^r {\cal O}_C) = \oplus^r \xi_{\cal G}$
is ample and has a section vanishing on $\mathbb P({\cal F})$.
\end{example}

\begin{corollary}\label{promin} Let $X$, ${\cal E}$ and $Z$ be as
in (\ref{setup}) with $dim Z \geq 2$. We assume that $Z$ is a
$\Proj$-bundle over a smooth
variety $W$ and also that $W$ is minimal.
Then $X$ is a $\Proj$-bundle over $W$ and ${\cal E}_{|F} = \oplus^r{\cal
O}_{\mathbb P}(1)$
for every fiber $F$ of $\phi : X \ra W$.
\end{corollary}

{\bf Proof.} \quad The assumption on the tautological bundle is not
necessary in this case as noted in (\ref{nohyposp}).

\begin{remark}
In case $r=1$ the last corollary shows that \cite[Conjecture
5.5.1]{BS} is true if $b\ge 3$, $X$ is smooth and $B$ is minimal.
\end{remark}

\begin{corollary}
Let $X$, ${\cal E}$ and $Z$ be as in (\ref{setup}) with dim~$Z \geq
2$. We assume that $Z$ has a scroll contraction $\f: Z \ra W$ with
dim~$W \leq 1$ (or equivalently that $Z$ is a $\Proj$-bundle over a
smooth variety $W$ of dimension $\leq 1$). Then $X$ is a
$\Proj$-bundle over $W$ and ${\cal E}_{|F} = \oplus^r{\cal
O}_{\mathbb P}(1)$ for every fiber $F$ of $\phi : X \ra W$ except
possibly for $W = \Proj^1$ and $Z = \Proj(\oplus^{(n-r)}{\cal
O}_{\mathbb P^1}) = \Proj^1 \times
\Proj^{(n-r-1)}$
or $Z =\Proj(\oplus^{(n-r-1)}{\cal O}_{\mathbb P^1}\oplus {\cal
O}_{\mathbb P^1}(1))$.
\end{corollary}

{\bf Proof.} \quad The corollary will follow if we prove that under
the assumptions $Z$ has only one extremal ray.

If dim~$W = 0$ then $Z = \Proj^{(n-r)}$ and thus $Z$ has only one
extremal ray. If dim~$W = 1$ then $\rho (Z) = 2$, thus $Z$ has one
extremal ray or it is Fano. But if $Z$ is a Fano manifold then $0 =
h^1({\cal O}_Z) = g(W)$, thus $W =
\Proj^1$.
Therefore we can assume that $Z= \Proj ({\cal E})$ for a vector bundle
${\cal E}$ on $\Proj^1$ with $rank({\cal E}) = s = n-r$ and $0 \leq
c_1({\cal E}) \leq s-1$.
But since $-K_Z = s\xi +(2- c_1({\cal E}))H$, with $\xi$ the tautological
bundle and
$H$ the pull back of a point in $\Proj^1$, if $c_1({\cal E}) \geq 2$ then
$\xi$ and
thus ${\cal E}$ would be ample. This is in contradiction with $c_1({\cal
E}) \leq s-1$.
Thus $ 0 \leq c_1({\cal E}) \leq 1$ which gives our claim.

\begin{remark}\label{notmin}  Even if $W$ is not minimal the $\Proj$-bundle contraction
can be sometime lifted to $X$. For instance if $Z$ is not Fano and
$Pic(W)= \Z$ the same proof of the above corollary gives the
lifting. In particular if $W = \Proj^m$ and ${\cal F}$ is a vector
bundle with $0 \leq c_1 {\cal F} \leq (m-1)$ which is not spanned
by global section then by \cite{SW1}, $Z$ is not Fano and the
$\Proj$-bundle contraction $Z=\Proj({\cal F}) \ra W$ can be lifted.
\end{remark}

\begin{remark} The above theorem and corollaries extend the results of Lanteri and Maeda;
their papers inspired and motivated ours. Summarizing we have a precise description of $X$
when  $Z$ is a scroll or a $\mathbb P$-bundle satisfying assumptions $(1.1)$ in the introduction (Theorem (\ref{mainscroll}) and corollary
(\ref{mainbundle})). We have a good result also if we drop the assumption but $W$ is minimal.
If $W$ is not minimal the situation is much more complicated; cases that occur if $\textrm{dim~}W \le 1$ are described in
the rest of this section; other simple cases are described in remark (\ref{notmin}).
\end{remark}

\begin{proposition}\label{excep1}
Suppose that $Z=\Proj(\oplus^{(n-r)}{\cal O}_{\mathbb P^1})=
\mathbb P^1 \times \mathbb P^{n-r-1}$ is the zero set of a section of an ample vector bundle
${\cal E}$ of rank r on a smooth $X$ of dimension $n$, then

\begin{enumerate}
\item $X$ is a $\mathbb P^{n-1}$-bundle over $\mathbb P^1$ and ${\cal E}_F = \oplus^r{\cal O}_F(1)$ for any fiber.
\item $X$ is a $\mathbb P^{r+1}$-bundle  over $\mathbb P^{n-r-1}$ and ${\cal E}_F = \oplus^r{\cal O}_F(1)$, for any fiber.
\item $(X,{\cal E})$ is $(\mathbb P^n,
\oplus^{n-3} {\cal O}_{\mathbb P}(1)\oplus{\cal O}_{\mathbb P}(2))$ or $(\mathbb Q^n,\oplus^{n-2}{\cal O}_{\mathbb Q}(1))$
and $Z$ is a smooth quadric surface.
\end{enumerate}
\end{proposition}

{\bf Proof.} \quad Let $Z = \Proj^1 \times \Proj^{(n-r-1)}$ and $p_1,p_2$ the
projections on the two factors; $Pic(Z) \simeq \mathbb Z
<p_1^*{\cal O}_{\mathbb P^1}(1)>\oplus~ \mathbb Z<p_2^*{\cal
O}_{\mathbb P^{n-r-1}}(1)>
 =: \mathbb Z <L_1> \oplus~ \mathbb Z <L_2>$; $p_1$ and $p_2$ are
two Fano-Mori contractions with supporting divisors $(aL_1,0)$ and
$(0,bL_2)$ respectively $(a,b >0)$; on $Z$ there is a third
Fano-Mori contraction, $p$, which is the contraction of $Z$ to a
point.

Suppose that $Z$ is the zero locus of a section of an ample vector
bundle $\cal E$ on $X$; by remark (\ref{nocor}) at least one of the
extremal contractions of $Z$ lifts to $X$.\par
\vspace{0.5 cm}
Suppose that $p_1$ lifts; as in the proofs of theorem (\ref{mainscroll}) and
corollary (\ref{mainbundle}) we obtain that $X$ is a $\mathbb
P^{n-1}$-bundle over $\mathbb P^1$ and ${\cal E}_F = \oplus^{r}
{\cal O}_F(1)$.\par
\vspace{0.5 cm}
Suppose now that $p_2$ lifts; as above we get that $X$ is a
$\mathbb P^{r+1}$-bundle over $\mathbb P^{n-r-1}$ and ${\cal E}_F
= \oplus^{r}{\cal O}_F(1)$.\par
\vspace{0.5 cm}

Finally, suppose that $p$ lifts. We start with the case in which $Pic(i):Pic(X)
\to Pic(Z)$ is an isomorphism and $\rho(X)=\rho(Z)$.
We want to show that in this
case both $p_1$ and $p_2$ would lift, but this is not possible, in view of theorem (\ref{manyrays}).\par
$X$ is a Fano variety and $-K_X = det \E + H$; this implies that the two extremal rays $R_1, R_2$ of $X$ have length $\ge r+1$;
If $R_i \quad (i= 1,2)$ is not nef, then the fibers of the associated contraction have dimension $\ge r+1$, and so, by proposition (\ref{intersezione}),
there exists a curve in the ray belonging to $Z$.
If $R_i \quad (i= 1,2)$ is nef, then the fibers of the associated contraction $P_i$ have dimension $\ge r$ by (\ref{iowi});  if a fiber has dimension $\ge r+1$
then its intersection with $Z$ contains at least a curve and we are done, so we can suppose that $P_i$ is equidimensional and all the fibers
have dimension $r$. Recalling that, for the general fiber 
$$ K_F = (K_X)_F = -(det \E)_F - H_F $$
we have that $(F,H_F) \simeq (\mathbb P^r, \Ol_{\mathbb P^r})(1))$, so, by theorem (\ref{pbundle}), $X$ is a $\mathbb P^r$-bundle on a normal
variety $W$.
This fact and the previous expression of $K_F$ yield  
$(F,\E_F) \simeq (\mathbb P^r, \oplus^{r}\Ol_{\mathbb P^r}(1))$ for every fiber; this implies that $Z \cap F$ is a point or contains a curve;
in the second case we are done, while the first is impossible, since $Z$ would be isomorphic to $W$, which has different
Picard number.\\

So suppose that $Pic(i):Pic(X) \to Pic(Z)$ is not an isomorphism or
$\rho(X) \ne \rho(Z)$; by (\ref{lef}) this is possible only for
dim$Z
=2$ and in this case we have $Pic(X) \simeq \mathbb Z$. By the
proof of theorem (\ref{lifting}) $K_X +det{\cal E} +2H = {\cal
O}_X$ for some ample line bundle $H$ such that $2H_Z = -K_Z$.
Note that there are curves $C$ on $Z$ such that $H.C = 1$, so $H$ is the ample
generator of $Pic(X)$; write $det {\cal E}= sH$; since the index of a Fano manifold is at most
$n+1$ we must have $s=r,r+1$. So $(X,\cal E)$ is either
$(\mathbb P^n,\oplus^{n-3} {\cal O}_{\mathbb P}(1)\oplus{\cal
O}_{\mathbb P}(2))$
 or $(\mathbb Q^n,\oplus^{n-2}{\cal O}_{\mathbb Q}(1))$
and $Z$ is a smooth quadric surface.\qquad $\square$\par

\begin{example}\label{ex1}

The effectiveness of case $3.$ is clear; to see the
effectiveness of case $1.$ consider the sequence
$$0 \lra \oplus^n {\cal O}_{\mathbb P^1} \lra \oplus^n(\op (a) \oplus \op(s-a)) \lra \oplus^n \op(s) \lra 0 $$
which is exact in view of \cite[Remark 1, p.170]{Ba2} and choose $a,s$ in such a way that
$ 0< a-s < a$; the construction in \cite[B.5.6]{Ful} applies and gives $\mathbb P^1 \times \mathbb P^{n-1}$ as a section
of the ample vector bundle ${\cal E} = \oplus^n{\xi_{\cal G}}$ on $X = \mathbb P({\cal G})$.
By the discussion above, it is clear
that in this example, the contraction $p_2$ cannot be lifted; $p_2$
is supported by $K_Z + H_Z= bL_2 \quad (b>0)$; recalling that $K_Z
=-2L_1 -(n-r)L_2$ we have that $H_Z= 2L_1 +(n-r+b)L_2$ is an ample
line bundle on $Z$ which cannot be the restriction on an ample line
bundle on $X$.
\end{example}

\begin{remark}
The effectiveness of case $2.$ looks uncertain; we note that for
$r=1$, i.e. in the case of ample divisors this is not possible by a
result of Sommese \cite[Theorem 5.2.1]{BS}.
\end{remark}

\begin{proposition}\label{F1}
Suppose that $Z=\Proj(\oplus^{(n-r-1)}{\cal O}_{\mathbb
P^1}\oplus {\cal O}_{\mathbb P^1}(1))$ is the zero set of a section
of an ample vector bundle ${\cal E}$ of rank r on a smooth $X$ of
dimension $n$, then

\begin{enumerate}
\item $X$ is a $\mathbb P^{n-1}$-bundle over $\mathbb P^1$ and ${\cal E}_F = \oplus^{n-1}{\cal O}_F(1)$ for any fiber.
\item $X$ is a scroll over $\mathbb P^{n-r}$ and ${\cal E}_F = \oplus^r{\cal O}_F(1)$ for the general fiber.
\item $X$ is both as in $1$ and in $2$ and it is $\mathbb P^1 \times
\mathbb P^{n-1}$, $r=1$ and ${\cal E}= {\cal O}(1,1)$.
\end{enumerate}
\end{proposition}

{\bf Proof.} \quad Let $Z =\Proj(\oplus^{(n-r-1)}{\cal O}_{\mathbb P^1}\oplus {\cal
O}_{\mathbb P^1}(1)) = $ blow-up of $\Proj^{(n-r)}$ along a linear
subspace of codimension $2$. The Mori cone of $Z$ is two
dimensional, and it is spanned by two extremal rays, $R_1$ and $R_2$;
let $p_i
\quad i=1,2$ be the extremal contraction associated to $R_i$; $p_1$
is the $\mathbb P^{n-r-1}$-bundle map on $\mathbb P^1$, while $p_2$
is the blow down to $\mathbb P^{n-r}$; moreover, let $p$ be the
contraction of $Z$ to a point, associated to the extremal face
generated by $R_1$ and $R_2$. By  (\ref{nocor}), there
is an extremal contraction that lifts to $X$.\par
\vspace{0.5 cm}
Suppose that $p_1$ lifts; as in the proof of (\ref{mainscroll})
and corollary (\ref{mainbundle}) we obtain that $X$ is a $\mathbb
P^{n-1}$-bundle over $\mathbb P^1$ and ${\cal E}_F = \oplus^r {\cal
O}_{\mathbb P^{n-1}}(1)$.\par
\vspace{0.5 cm}
Suppose now that $p_2$ lifts to an extremal contraction $P_2$ on $X$.

\begin{lemma}
The contraction $P_2$ is of fiber type.
\end{lemma}

{\bf Proof.} \quad Assume by contradiction that $P_2$ is birational. \par
{\bf Claim.} \quad $P_2$ is divisorial and it is the blow up of a smooth $X'$ along a subvariety of codimension $r+2$.\par

The contraction $P_2$ is supported by a divisor of the form $K_X + det \E + H$, and this implies dim~$F \ge l(R) \ge r+1$,
so, reasoning as in the proof of (\ref{agreement}), we can prove that ${P_2}_Z$ has connected fibers.
The fact that ${P_2}_Z$ factors through $p_2$ and the normalization of 
${P_2}_Z(Z)$ yields that all the nontrivial fibers of ${P_2}_Z$ have dimension 1; this implies that all the nontrivial
fibers of $P_2$ have dimension at most $r+1$ by proposition (\ref{intersezione}); 
so we can conclude that the contraction $P_2$ is equidimensional and all the nontrivial
fibers have dimension $r+1$; by (\ref{iowi}) we
have
$$dim E(P_2) \ge l(R) - dim F + n - 1 \ge n-1$$
so $P_2$ is divisorial and dim $P_2(E(P_2))= (n-1)-(r+1)= n-r-2$.
In view of $(\ref{lef})$ the contraction $P_2$ is elementary, and so we can
apply \cite[Theorem 4.1 and Corollary 4.11]{AWD}, proving the
claim.\par
The Picard group of $Z$ is generated by the tautological line bundle $\xi$ and $f$, a fiber of
the projection on $\mathbb P^1$, but also by $E(p_2)$ and
$p_2^*{\cal O}(1)$; we have $ p_2^*{\cal O}(1) = \xi $ and $ E(p_2)  = \xi - f$.
Using the isomorphism $Pic(i):Pic(X) \to Pic(Z)$ and the fact that
$p_2$ and $P_2$ are elementary contractions, we get that the
restriction omomorphism $Pic(j):Pic(X') \to Pic(\mathbb P^{n-r})$
is an isomorphism. So there exists an ample generator of $Pic(X')$,
$H'$ whose restriction to $\mathbb P^{n-r}$ is ${\cal O}_{\mathbb
P}(1)$. $Pic(X)$ is thus generated by $P_2^*H$ and $E(P_2)$. Write
$K_{X'} = k H'$ for some $k \in \mathbb Z$. Let $l$ be a line in
$f$; we obtain

$$ -(n-r) = K_Z.l = (K_X +det {\cal E})_Z.l =$$
$$= (P_2^*K_{X'}+ (r+1)E(P_2) + det {\cal E})_Z.l = $$
$$ = (kP_2^*H' + (r+1)E(P_2) + det {\cal E})_Z.l = $$
$$= ((k+r+1)\xi -(r+1)f+det{\cal E}_Z).l =$$
$$ k+r+1 +det {\cal E}_Z.l \ge k+2r+1$$

and so
$$k \le -n-r-1$$

which is absurd, since the index of a Fano variety is not greater
than $n+1$, and the lemma is proven.\par \vspace{0.5 cm}
On $X$ we thus have a fiber type elementary contraction, supported
by an high multiple of $K_X + det {\cal E} +H$ on a variety of
dimension $n-r$ by (\ref{surjectivity}). For the general fiber of
$P_2$ we have

$$ K_F + det{\cal E}_F+H_F = {\cal O}_F$$

and so $(F,{\cal E}_F) = (\mathbb P^r, \oplus^r {\cal O}_{\mathbb
P}(1))$; therefore $Z \cap F$ is a point for the generic fiber, and
thus ${P_2}_Z$ is generically one-to-one and therefore coincides with $p_2$ (see
(\ref{agreement})). So the conclusion is that, if $p_2$ lifts, $X$
is a scroll over $\mathbb P^{n-r}$.\par\vspace{0.5 cm}

As a final case suppose that $p$ lifts; if $Pic(i):Pic(X) \to
Pic(Z)$ is an isomorphism and $\rho(X) = \rho (Z)$; as in the proof of proposition (\ref{excep1})
we can prove that also $p_1$ and $p_2$ lift, so
$X$ is a Fano variety which has a $\mathbb P^{n-1}$-bundle
contraction on $\mathbb P^1$ and a scroll contraction on $\mathbb
P^{n-r}$. The only possibility is that $r=1$ and $X=\mathbb
P_{\mathbb P^1}(\oplus^r{\cal O}) =
\mathbb P^1 \times \mathbb P^{n-1}$.\par

If $Pic(i)$ is not an isomorphism or $\rho(X) \ne \rho(Z)$, by theorem (\ref{lef}) $Pic(X)
\simeq \mathbb Z$ and dim$Z = 2$, so $Z = \mathbb F_1$ and $X$ is a Fano variety; this case is ruled out in
\cite[Section 2]{LM3}.\qquad $\square$

\begin{example}\label{ex2}
Cases $1.$ and $3.$ are effective; examples for the first case can
be constructed as in example (\ref{ex1}) starting with the exact
sequence $$ 0 \lra \oplus^n\op \lra \oplus^{n-1}(\op(a)\oplus \op(s-a))\oplus \op(a) \oplus \op(s+1-a)
\lra \dots\hfill$$
$$\hfill \dots \lra \oplus^{n-1}\op(s) \oplus \op(s+1) \lra 0$$
while in case 3. easy computations show that a smooth $Z$ in the linear system $|{\cal O}(1,1)|$
must be a Fano variety with a $\mathbb P$-bundle contraction on $\mathbb P^1$ and a birational contraction on
$\mathbb P^{n-1}$.\par
\end{example}

\begin{remark} The example of case 3. shows that the assumption $\tau \ge 2$ in (\ref{birubiru}) is necessary. \end{remark}

\begin{remark}
As in the case  $Z=\mathbb P^1 \times \mathbb P^{n-r-1}$
the effectiveness of case $2.$ looks uncertain; we note that for
$r=1$, i.e. in the case of ample divisors this is not possible except for the trivial case $X= \mathbb P^{n-1} \times \mathbb P^1$,
${\cal E} = \Ol_{X= \mathbb P^{n-1} \times \mathbb P^1}(1,1)$.
This is a very well known fact that descends by a result of Sommese \cite[Theorem 5.2.1]{BS} and classical results of B\v{a}descu,
\cite{Ba2}, \cite{Ba3}.\end{remark}

\section{Quadric fibrations and quadric bundles}

\begin{theorem}\label{mainquadric}
Let $X$, $\cal E$ and $Z$ be as in (\ref{setup}) with $dim Z \geq
3$. We assume that $Z$ has a quadric fibration contraction $\f: Z
\ra W$ with respect to an ample line bundle on $Z$, $H_Z$, which is
the restriction of an ample line bundle $H$ on $X$. Then $X$ has a
Fano-Mori contraction $\phi:X\to W$ which is of fiber type and with
supporting divisor $D=K_X+det {\cal E}+(n-m-r)H$ with $n =$ dim~
$X$ and $m=$ dim~$W$. For the general fiber of $\phi$, $F$ we have
either $(F,{\cal E}_F)
\simeq (\mathbb P^{n-m},\oplus^{r-1}{\cal O}_{\mathbb
P}(1)\oplus{\cal O}_{\mathbb P}(2))$ or $(F,{\cal E}_F) \simeq
(\mathbb Q^{n-m},\oplus^r {\cal O}_{\mathbb Q}(1))$.
\\
If $\varphi$ is elementary then also $\phi$ is elementary and it is either a
scroll contraction or a quadric fibration contraction (i.e. it is
supported by the divisor $K_X+(n-m+1)H$ or by the divisor $K_X
+(n-m)H$).
The last result holds replacing the assumption on $\varphi$ with the strongest assumption
dim~$X \geq 2 m+1 = 2 dim W +1$ (this is always the case if $dim W \ge 1$).\par
\end{theorem}

{\bf Proof.} \quad The morphism $\varphi$ is a contraction
supported by $K_Z +(n-r-m)H_Z$, so, applying theorem
(\ref{lifting}), we get a contraction $\phi:X\to Y$, defined by an
high multiple of $D=K_X+det {\cal E}+ (n-m-r)H$; this contraction
is of fiber type and $Y = W$ by proposition (\ref{agreement}). Let
$F$ be a general fiber of $\phi$; then $F$ is a smooth Fano
manifold of dimension $n-m$ such that $-K_F=({det \cal
E}+(n-m-r)H)_F$. Thus, by \cite[Theorem 0.1]{PSW} applied to the
ample vector bundle ${\cal E}_1 = ({\cal E}
\oplus^{n-m-r}H)_F$ either $F=\mathbb P^{n-m}$ and $\cal E$
restricted to it is $\oplus^{r-1}{\cal O}_{\mathbb P}(1)\oplus
{\cal O}_{\mathbb P}(2)$ or $F=\mathbb Q^{n-m}$ and $\cal E$
restricted to it is $\oplus^r{\cal O}_{\mathbb Q}(1)$. Moreover,
for any line in a general fiber
\begin{equation}\label{eqqb} (det {\cal E} -(r+\epsilon)H).l = 0 \end{equation}
with $\epsilon = 1,0$. By theorem (\ref{lef}) $\rho(X/W)=\rho(Z/W)$
and so, if $\f$ is elementary, also $\phi$ is so and, by
(\ref{eqqb}) $det {\cal E}= (r+\epsilon)H +\phi^*B$, that is $\phi$
is supported by $K_X +(n-m+\epsilon)H$. Assume now that dim~$X \ge
2$ dim~$W +1$; by proposition (\ref{elementary}) the contraction
$\phi:X\to W$ is an elementary contraction and so, again by
(\ref{lef}) also $\varphi:Z \to W$ had to be elementary. $\square$\par

\begin{corollary}
Assume now that there exists a vector bundle ${\cal F}$ on $W$ and
an embedding of $Z$ into $\mathbb P({\cal F})$ as a divisor of
relative degree $2$; assume moreover that $K_Z +(n-r-m)(\xi_{\cal
F})_Z$ is a good supporting divisor of a quadric bundle elementary
contraction and $(\xi_{\cal F})_Z$ is the restriction of an ample
line bundle on $X$. Then either there exists an ample vector bundle
${\cal G}$ of rank $n-m+1$ such that $X=\mathbb P({\cal G})$ or
there exists a vector bundle ${\cal G}$ of rank $n-m+2$ and an
embedding of $X$ into $\mathbb P ({\cal G})$ as a divisor of
relative degree two.
\end{corollary}

{\bf Proof.} If $Z$ is a quadric bundle then $\phi$ is
equidimensional; in fact if it has any fiber of dimension $> n-m$
then, by proposition (\ref{intersezione}), even $Z \ra W$ should have
a fiber of dimension $>(n-m-r)$. Since $\f: Z \ra W$ is elementary
$\phi$ is a scroll or a quadric fibration with the respect to $H$.
We conclude by (\ref{pbundle}) and (\ref{qbundle}).\par

\begin{remark} Theorem (\ref{mainquadric}) extends \cite[Theorem C]{LM2} and rules out the doubtful case $(3)$.\end{remark}

\section{Del Pezzo fibrations}

\begin{theorem}\label{maindelpezzo}

Let $X$, $\cal E$ and $Z$ be as in (\ref{setup}) with $dim Z \geq
3$. We assume that $Z$ has a del Pezzo fibration contraction $\f: Z
\ra W$ with respect to an ample line bundle on $Z$, $H_Z$, which is
the restriction of an ample line bundle $H$ on $X$. Then $X$ has a
Fano-Mori contraction $\phi:X\to W$ which is of fiber type and with
supporting divisor $D=K_X+det {\cal E}+(n-m-r-1)H$, with $n =$ dim~
$X$ and $m=$ dim~$W$. If $n-m\ge 5$, for the general fiber of
$\phi$, $F$, we have either $(F,{\cal E}_F)
\simeq (\mathbb P^{n-m},\oplus^{r-1}{\cal O}_{\mathbb
P}(1)\oplus{\cal O}_{\mathbb P}(3))$, $(F,{\cal E}_F)
\simeq (\mathbb P^{n-m},\oplus^{r-2}{\cal O}_{\mathbb
P}(1)\oplus^2{\cal O}_{\mathbb P}(2))$, or $(F,{\cal E}_F) \simeq
(\mathbb Q^{n-m},\oplus^{r-1} {\cal O}_{\mathbb Q}(1)\oplus{\cal
O}_{\mathbb Q}(2))$ or $F$ is a del Pezzo manifold with $b_2=1$ and
${\cal E}_F \simeq \oplus^r{\cal O}_F(1)$, where ${\cal O}_F(1)$ is
the ample generator of $Pic(F)$. \\
If $\varphi$ is elementary then $\phi$ is elementary and it is either a scroll contraction or a
quadric fibration contraction or a del Pezzo fibration contraction
(i.e. it is supported by the divisor $K_X+(n-m+\epsilon)H$ with $\epsilon = 1,0$ or $-1$).
The last result holds replacing the assumption on $\varphi$ with the strongest assumption
 $dim X = n \geq 2 m+3 = 2 dim W +3$.\par
\end{theorem}

{\bf Proof} \quad The morphism $\varphi$ is a contraction supported
by $K_Z +(n-r-m-1)H_Z$, so, applying theorem (\ref{lifting}),
we get a contraction $\phi:X\to Y$, defined by an high multiple of
$D=K_X+det {\cal E}+ (n-m-r-1)H$; this contraction is of fiber type
and $Y = W$ by proposition (\ref{agreement}). Let $F$ be a general
fiber of $\phi$; then $F$ is a smooth Fano manifold of dimension
$n-m$ such that $-K_F=({det \cal E}+(n-m-r)H)_F$. Thus, applying
\cite[Main theorem]{PSW} to the ample vector bundle ${\cal E}_1
= ({\cal E}\oplus ^{n-m-r-1}H)_F$ we get the description of $F$ and ${\cal E}_F$.
Moreover, for any line in a general fiber
\begin{equation}\label{eqdp}
(det{\cal E}-(r+\epsilon)H).l = 0 \end{equation}
with $\epsilon = 1,0$ or $-1$.\par
By theorem (\ref{lef})
$\rho(X/W)=\rho(Z/W)$ and so, if $\f$ is elementary, also $\phi$ is
so, by (\ref{eqdp}) $det {\cal E}= (r+\epsilon)H +\phi^*B$, that is
$\phi$ is supported either by $K_X +(n-m+\epsilon)H$. Assume now
that dim~$X \ge 2$ dim~$W +3$; by proposition (\ref{elementary})
the contraction $\phi:X\to W$ is an elementary contraction and so,
again by (\ref{lef}) also $\varphi:Z
\to W$ had to be elementary. $\square$\par

\begin{remark}\label{dp4}
If $n-m=4$ the rank of ${\cal E}$ can be $1$ or $2$ and, according
to \cite[proposition 7.4]{PSW} the possibilities for the general fiber are those listed
in the theorem plus
\begin{enumerate}
\item $(F,{\cal E}_F) \simeq (\mathbb P^2 \times \mathbb P^2,\oplus^2 {\cal O}(1,1)).$
\item $(F,{\cal E}_F) \simeq (\mathbb P^2 \times \mathbb P^2,{\cal O}(1,1)).$
\item $(F,{\cal E}_F) \simeq (\mathbb Q^4, S(2)).$
\item $F$ is a Fano 4-fold with $b_2 = 1$ and index $1$.
\end{enumerate}
\end{remark}
\begin{remark}\label{dp3}
If $n-m=3$ the rank of ${\cal E}$ must be $1$ and, according to
\cite[Theorem 0.4]{PSW} the possibilities for the general fiber are those listed in the
theorem plus
\begin{enumerate}
\item $(F,{\cal E}_F) \simeq (\mathbb P^2 \times \mathbb P^1,{\cal O}(2,1)).$
\item $(F,{\cal E}_F) \simeq (\mathbb P^1 \times \mathbb P^1 \times \mathbb P^1,{\cal O}(1,1,1)).$
\item $F$ is a del Pezzo manifold with $b_2\ge 2$ and ${\cal E}={\cal O}_F(1).$
\end{enumerate}
\end{remark}
\begin{corollary}
Let $X$, ${\cal E}$ and $Z$ be as in (\ref{setup}) and let $Z$ be a
del Pezzo manifold with $b_2=1$. then one of the following occurs
\begin{enumerate}
\item $X\simeq \mathbb P^n$ and ${\cal E}$ is either $\oplus^2\opn(2)\oplus^{r-2}\opn(1)$
or $\opn(3)\oplus^{r-1}\opn(1)$.
\item $X \simeq \mathbb Q^n$ and ${\cal E}$ is $\oqn(2)\oplus^{r-1} \oqn(1)$.
\item $X$ is a del Pezzo manifold with $b_2=1$ and ${\cal E} \simeq \oplus^r{\cal O}_X(1)$
where ${\cal O}_X(1)$ is the ample generator of $Pic(X)$.
\end{enumerate}
\end{corollary}
{\bf Proof.} \quad The hypothesis on $H$ is not necessary in this
case as noted in (\ref{nohypo}) and the cases in (\ref{dp4}) and
(\ref{dp3}) are ruled out, because of the isomorphism $Pic(Z) \simeq Pic(X)
\simeq \mathbb Z$.

\section{Some final considerations}

Using the same arguments we can consider the case in which $Z$ has
an extremal contraction on $W$ whose general fiber $F$ is a Fano
variety of index $\le$ dim$F -2$. However, in these cases it is very
difficult to provide a good description of the vector bundle $\cal
E$ and to construct non trivial examples. These difficulties
already show up in the case in which $W$ is a point and $Pic(Z)
\simeq Z$; we recall that a line on $Z$ in this case is a rational
curve which is a line with respect to a the generator of $Pic(Z)$.
The existence of a line on $Z$ is proved if the index of $Z$ is
$\geq (n-2)$, by recent results of M. Mella. A line exists also if
$2 \textrm{~index~}(Z) > \textrm{dim~}Z +1$. The following
proposition summarize the simplest cases.

\begin{proposition}
Let $X,{\cal E}$ and $Z$ be as in (\ref{setup}); we assume that $Z$
is a Fano variety of dimension $\ge 2$ with $Pic(Z) \simeq
\mathbb Z$ and that $Z$ has a line. Then $X$ is a Fano variety
with $Pic(X) \simeq \mathbb Z$ and coindex$(Z) \ge$ coindex $(X)$.
\end{proposition}

{\bf Proof.} Let $H$ be a generator of $Pic(X) = Pic(Z)$ and let
$s$ and $\tau$ be positive integers such that $det{\cal E} = s H$
and $-K_Z = \tau H_Z$. Therefore, by adjunction formula and
(\ref{agreement}), we have that
$$K_X + det {\cal E}+\tau H =K_X+(s+\tau)H={\cal O}_X,$$
thus $X$ is a Fano manifold. If $C$ is a line of $Z$ then $det{\cal
E} ^.C = s H^.C = s$; thus, by (\ref{det}), $s \geq r$. In
particular this gives $n+1 \geq \textrm{index} (X) = s+\tau \geq r+
\tau$.

\medskip
\scriptsize{DIPARTIMENTO DI MATEMATICA, UNIVERSIT\'A DI TRENTO, 38050 POVO (TN), ITALIA \\
{\it E-mail address: andreatt@science.unitn.it} \\
\null \\
DIPARTIMENTO DI MATEMATICA, UNIVERSIT\'A DI TRENTO, 38050 POVO (TN), ITALIA \\
{\it E-mail address: occhetta@science.unitn.it}} \\


\begin{thebibliography}{KMM92}

\bibitem[ABW91]{ABW1}
M.~Andreatta, E.~Ballico, and J.~A. Wi{\'s}niewski.
\newblock On contractions of smooth algebraic varieties.
\newblock {\em preprint UTM 344}, ~(1991).

\bibitem[ABW93]{ABW3}
\rule[.5ex]{2cm}{.2mm}.
\newblock Two theorems on elementary contractions.
\newblock {\em Math. Ann.}, {\bf 297}~(1993),~191--198.


\bibitem[AW93]{AWD}
M.~Andreatta and J.~A. Wi{\'s}niewski.
\newblock A note on nonvanishing and applications.
\newblock {\em Duke Math. J.}, {\bf 72}~(1993),~739--755.

\bibitem[AW97]{AW}
\rule[.5ex]{2cm}{.2mm}.
\newblock A view on contraction of higher dimensional varieties.
\newblock In {\em Algebraic Geometry -- Santa Cruz 1995}, volume~62 of {\em
 Proc. Sympos. Pure Math.},  AMS, Providence, RI,
 1997, 153--183.

\bibitem[B{\v a}d81]{Ba2}
L.~B{\v a}descu.
\newblock On ample divisors, II.
\newblock In {\em Week of Algebraic Geometry, Proc. Bucharest}, volume~40 of
  {\em Texte Math.},  Teubner, Leipzig, 1981, 12--32.

\bibitem[B{\v a}d82a]{Ba1}
\rule[.5ex]{2cm}{.2mm}.
\newblock On ample divisors.
\newblock {\em Nagoya Math J.}, {\bf 86}~(1982),~155--171.

\bibitem[B{\v a}d82b]{Ba3}
\rule[.5ex]{2cm}{.2mm}.
\newblock The projective plane blown up at a point as an ample divisor.
\newblock {\em Atti Accad. Ligure Sci. Lett.}, {\bf 38}~(1982),~88--92.
.
\bibitem[BS92]{BS}
M.~C. Beltrametti and A.~J. Sommese.
\newblock On the adjunction theoretic classification of polarized varieties.
\newblock {\em J. reine und angew. Math}, {\bf 427}~(1992),~157--192.

\bibitem[BS95]{BS2}
\rule[.5ex]{2cm}{.2mm}.
\newblock {\em The adjunction theory of complex projective varieties},
  volume~16 of {\em Exp. Math.}
\newblock de Gruyter, Berlin, 1995.

\bibitem[BSW90]{BSW}
M.~C. Beltrametti, A.~J. Sommese and J.~A. Wi{\'s}niewski.
\newblock Results on varieties with many lines and their application to adjunction theory.
\newblock In {\em Complex algebraic varieties, Bayreuth 1990},  {\em Lecture Notes
in  Math.}, {\bf 1507}~   Springer-Verlag,  1992, 16--38.

\bibitem[Fuj80]{Fu2}
T.~Fujita.
\newblock Hyperplane section principle.
\newblock {\em Journal of the Math. Soc. of Japan}, {\bf 32}~(1980),~153--169.

\bibitem[Fuj90]{Fub}
\rule[.5ex]{2cm}{.2mm}.
\newblock {\em Classification theory of polarized varieties}, volume 155 of
  {\em London Math. Soc. Lecture Notes Series}.
\newblock Cambridge Univ. Press., Cambridge, 1990.

\bibitem[Ful84]{Ful}
W.~Fulton.
\newblock {\em Intersection theory}, volume~2 of {\em Ergebnisse der Math.}
\newblock Springer Verlag, Berlin Heidelberg New York Tokio, 1984.

\bibitem[KMM87]{KMM}
Y.~Kawamata, K.~Matsuda, and K.~Matsuki.
\newblock Introduction to the minimal model program.
\newblock In {\em Algebraic geometry, Sendai}, volume~10 of {\em Adv. Studies
  in Pure Math.},  Kinokuniya-North-Holland, 1987, 283--360.

\bibitem[KMM92]{KoMiMo}
J.~Koll\'ar, Y.~Miyaoka, and S.~Mori.
\newblock Rational connectedness and boundedness of {Fano} manifolds.
\newblock {\em J. Diff. Geom.}, {\bf 36}~(1992),~765--779.


\bibitem[LM95]{LM1}
A.~Lanteri and H.~Maeda.
\newblock Ample vector bundles with sections vanishing on projective spaces or
  quadrics.
\newblock {\em Intern. J. Math.}, {\bf 6}~(1995),~587--600.

\bibitem[LM96]{LM2}
\rule[.5ex]{2cm}{.2mm}.
\newblock Ample vector bundles characterizations of projective bundles and
  quadric fibrations over curves.
\newblock In {\em Higher dimensional complex varieties (Trento, 1994)}, 
 de Gruyert, 1996,  247--259.

\bibitem[LM97]{LM3}
\rule[.5ex]{2cm}{.2mm}.
\newblock Geometrically ruled surfaces as zero loci of ample vector bundles.
\newblock {\em Forum Math.}, {\bf 9}~(1997),~1--15.

\bibitem[Mor79]{Mo2}
S.~Mori.
\newblock Projective manifolds with ample tangent bundle.
\newblock {\em Ann. Math.}, {\bf 110}~(1979),~593--606.

\bibitem[Mor82]{Mo1}
S.~Mori.
\newblock Threefolds whose canonical bundles are not numerically effective.
\newblock {\em Ann. Math.}, {\bf 116}~(1982),~133--176.

\bibitem[Pet90]{Pe1}
T.~Peternell.
\newblock A characterization of {$\mathbb P^n$} by vector bundles.
\newblock {\em Math. Zeitschrift}, {\bf 205}~(1990),~487--490.

\bibitem[PSW92]{PSW}
T.~Peternell, M.~Szurek, and J.~A. Wi{\'s}niewski.
\newblock Fano manifolds and vector bundles.
\newblock {\em Math. Ann.}, {\bf 294}~(1992),~151--165.

\bibitem[Som76]{So1}
A.~J. Sommese.
\newblock On manifolds that cannot be ample divisors.
\newblock {\em Math. Ann.}, {\bf 221}~(1976),~55--72.

\bibitem[Som78]{So2}
\rule[.5ex]{2cm}{.2mm}.
\newblock Submanifolds of abelian varieties.
\newblock {\em Math. Ann.}, {\bf 233}~(1978),~229--256.

\bibitem[Som86]{So}
\rule[.5ex]{2cm}{.2mm}.
\newblock On the adjunction theoretic structure of projective varieties.
\newblock In {\em Complex Analysis and Algebraic Geometry, Proceedings
  G{\"o}ttingen , ed. H.Grauert}, volume 1194 of {\em Lecture Notes in Math.},
 Springer-Verlag, 1986, 175--213. 

\bibitem[SW90a]{SW2}
M.~Szurek and J.~A. Wi{\'s}nieswki.
\newblock Fano bundles over {$\mathbb P^3$} and {$\mathbb Q^3$}.
\newblock {\em Pacific J. of Math.}, {\bf 141}~(1990),~197--208.

\bibitem[SW90b]{SW1}
\rule[.5ex]{2cm}{.2mm}.
\newblock On {Fano} manifolds which are {${\mathbb P^k}$}-bundles over
  {${\mathbb P^2}$}.
\newblock {\em Nagoya Math. J.}, {\bf 120}~(1990),~89--101.

\bibitem[Wi{\'s}90]{Wi2}
J.~A. Wi{\'s}niewski.
\newblock On a conjecture of Mukai.
\newblock {\em Manuscripta Math}, {\bf 68}~(1990),~135--141.

\bibitem[Wi{\'s}91]{Wi1}
J.~A. Wi{\'s}niewski.
\newblock On contractions of extremal rays of Fano manifolds.
\newblock {\em J. reine und angew. Math}, {\bf 417}~(1991),~141--157.



\end{thebibliography}
\end{document}